\documentclass[a4 paper,12pt]{article}
\usepackage{epsfig,graphics}
\usepackage[english]{babel}
\usepackage{amssymb}
\usepackage{amsmath}
\usepackage{array}
\newcolumntype{T}{>{\tiny}l}
\newcolumntype{H}{>{\Huge}l}
\newcolumntype{S}{>{\small}l}
\usepackage{multirow}
\usepackage{array}
\usepackage{geometry}
\usepackage{graphicx}

\begin{document}

\vspace{.2in}\parindent=0mm

\begin{flushleft}
 {\bf\Large {Biorthogonal Wavelets on the Spectrum }}

 \parindent=0mm \vspace{.4in}
  {\bf{Owais Ahmad$^{\star}$ and F. A. Shah}}
\end{flushleft}

\parindent=0mm \vspace{.1in}
{{\it\small$^{\star}${Department of Mathematics,National Institute of Technology, Srinagar-190006, Jammu and Kashmir, India.
E-mail: $\text{siawoahmad@gmail.com}$}}

\parindent=0mm \vspace{.1in}
 {{\it\small{Department of  Mathematics,  University of Kashmir, South Campus, Anantnag-192 101, Jammu and Kashmir, India. }}

\begin{quote}
  {\small Abstract: A generalization of Mallat's classic theory of multiresolution analysis based on the theory of spectral pairs was considered by Gabardo and Nashed (J. Funct. Anal. 158,  209-241, 1998). In this article, we introduce the notion of biorthgonoal nonuniform multiresolution analysis on the spectrum $\Lambda=\left\{0, r/N\right\}+2\mathbb Z$, where $N\ge 1$ is an integer and $r$ is an odd integer with $1\le r\le 2N-1$ such that $r$ and $N$ are relatively prime. We first establish the necessary and sufficient conditions for the translates of a single function to form the Riesz bases for their closed linear span. We provide the complete characterization for the biorthogonality of  the translates of scaling functions of two nonuniform multiresolution analysis and the associated biorthogonal wavelet families. Furthermore, under the mild assumptions on the scaling functions and the corresponding wavelets associated with nonuniform multiresolution analysis, we show that the wavelets can generate Reisz bases.}
\end{quote}

\parindent=0mm \vspace{.0in}

\parindent=0mm \vspace{.1in}

\parindent=0mm \vspace{.1in}
{\bf{1. Introduction}}

\parindent=0mm \vspace{.1in}
The theory of continuous and discrete wavelet transforms have emerged as a widely used tool  in various disciplines of science and engineering including    image processing, spectrometry,  turbulence, computer graphics, optics and electromagnetism, telecommunications, DNA sequence analysis,  quantum physics,   statistics, solution of partial and ordinary differential equations. Orthogonality has long been assumed as a key property in virtually all standard approaches when analyzing or synthesising signals. A higher-level signal processing technique involves the concept of biorthogonality in which two (cross-orthogonal) sets are used: one for the analysis and the other one synthesis. During the late 1990's, biorthogonal wavelets brought a major breakthrough into image compression, thanks to their natural feature of concentrating energy in a few transform coefficients. In traditional wavelet theory, biorthogonal wavelets have many advantages over orthogonal wavelets, by relaxing orthonormal to biorthogonal, additional degrees of freedom are added to design problems. Biorthogonal wavelets in $L^2(\mathbb R)$ were investigated by Bownik and Garrigos \cite{BG}, Cohen et al.\cite{CDF}, Chui and Wang \cite{CW} and many  others.

\parindent=8mm \vspace{.1in}

Multiresolution analysis is an important mathematical tool since it provides a natural framework for understanding and constructing discrete wavelet systems. The concept of MRA has been extended in various ways in recent years. These concepts are generalized to  $L^2\big(\mathbb R^d\big)$, to lattices different from  $\mathbb Z^d$, allowing the subspaces of MRA to be generated by Riesz basis instead of orthonormal basis, admitting a finite number of scaling functions, replacing the dilation factor 2 by an integer $M\geq 2$ or by an expansive matrix $A\in GL_{d}(\mathbb R)$ as long as $A\subset A\mathbb Z^d$. All these concepts are developed on regular lattices, that is the translation set is always a group. Recently, Gabardo and Nashed \cite{GN, GY} considered a generalization of Mallat's  celebrated theory of  MRA based on spectral pairs, in which the translation set acting on the scaling function associated with the  MRA to generate the subspace $V_0$ is no longer a group, but is the union of $\mathbb Z$ and a translate of $\mathbb Z$. They call it \textit{nonuniform multiresolution analysis} (NUMRA). A parallel development on local fields of positive characteristic was carried by  Shah and his colleagues \cite{SA,SB}. The concept of biorthogonal wavelets associated with MRA and NUMRA were studied by various researchers \cite{BG,LC,SM}.

\pagestyle{myheadings}

\parindent=8mm \vspace{.1in}
In this paper, we introduce the notion of biorthogonal wavelets on the spectrum and  obtain the necessary and sufficient conditions for the translates of a single function to form the Riesz bases for their closed linear span. We also provide a complete characterization for the biorthogonality of  the translates of scaling functions of two NUMRA's and the associated biorthogonal wavelet families. Moreover, under  mild assumptions on the scaling functions and the corresponding wavelets,  we show that the nonuniform wavelets can generate Reisz bases for $L^2(\mathbb R)$.

\parindent=8mm \vspace{.1in}
The article is structured in the following manner. In Section 2, we recall the basic definitions of MRA and NUMRA . In Section 3, we establish necessary and sufficient conditions for the translates of a function to form a Riesz basis for its closed linear span. In the concluding Section, we show that the wavelets associated with dual MRA's are biorthogonal and
generate Riesz bases for $L^2(\mathbb R)$.

\parindent=0mm \vspace{.2in}
{\bf{2. Preliminaries  }}

\parindent=0mm \vspace{.1in}
{\bf{Definition 2.1.}}  An MRA of $L^2({\mathbb R})$ is a  sequence of closed subspaces $\{V_j:j\in \mathbb Z\}$ of $L^2({\mathbb R})$ satisfying the following properties:

\parindent=0mm \vspace{.1in}
(a)\quad $V_j \subset V_{j+1}\; \text{for all}\; j \in \mathbb Z;$

\parindent=0mm \vspace{.1in}
(b)\quad $\bigcup_{j\in \mathbb Z}V_j\;\text{is dense in}\;L^2({\mathbb R});$

\parindent=0mm \vspace{.1in}
(c)\quad $\bigcap_{j\in \mathbb Z}V_j=\{0\};$

\parindent=0mm \vspace{.1in}
(d)\quad $f(x) \in V_j\; \text{if and only if}\;f(2x) \in V_{j+1}\; \text{for all}\; j \in \mathbb Z;$

\parindent=0mm \vspace{.1in}
(e) There is a function $\phi \in V_0$, called the {\it{scaling function}}, such that $\left\{\phi(x-k): k\in \mathbb Z\right\}$ forms an orthonormal basis for $V_0$.

\parindent=8mm \vspace{.1in}
According to the standard scheme for construction of MRA-based wavelets, for each $j$, we define a space $W_{j}$ ({\it wavelet space}) as the orthogonal complement of $V_{j}$ in $V_{j+1}$, i.e., $V_{j+1}=V_{j}\oplus W_{j}, \, j\in\mathbb Z$, where $W_{j}\perp V_{j}, \, j\in\mathbb Z$. It is not difficult to see that
$$f(x)\in W_{j} \quad \text{if and only if}\quad f(2x)\in W_{j+1},\quad j\in\mathbb Z.\eqno(2.1)$$

\parindent=0mm \vspace{.1in}
Moreover, they are mutually orthogonal, and we have the following orthogonal decompositions:
$$L^2({\mathbb R})= \bigoplus_{j\in\mathbb Z} W_{j}=V_{0}\oplus \left(\bigoplus_{j\ge 0}W_{j}\right).\eqno(2.2)$$

\parindent=0mm \vspace{.1in}
For an integer $N\ge 1$ and an odd integer $r$ with $1 \le r \le 2N-1$ such that $r$ and $N$ are relatively prime, we define
$$\Lambda=\left\{0, \dfrac{r}{N}\right\}+2\mathbb Z=\left\{ \dfrac{rk}{N}+2n: n\in \mathbb Z, k=0,1\right\}.\eqno(2.3)$$

\parindent=0mm \vspace{.1in}
It is easy to verify that $\Lambda$ is not necessarily a group nor a uniform discrete set, but is the union of $\mathbb Z$ and a translate of $\mathbb Z$.
 Moreover, the set $\Lambda$ is the spectrum for the spectral set $\Gamma=\big[0, \frac{1}{2}\big)\cup \big[\frac{N}{2}, \frac{N+1}{2}\big)$ and the
pair $(\Lambda,\Gamma)$ is called a {\it spectral pair} \cite{GN}.

\pagestyle{myheadings}

\parindent=0mm \vspace{.1in}
{\bf{Definition 2.2.}} For an integer $N \ge 1$ and an odd integer $r$ with $1\leq r \leq 2N-1$ such that $r$ and $N$ are relatively prime, an associated  nonuniform MRA  is a sequence of closed subspaces $\left\{V_j: j\in\mathbb Z\right\}$ of $L^2({\mathbb R})$ such that the following properties hold:

\parindent=0mm \vspace{.1in}
(a)\quad $V_j \subset V_{j+1}\; \text{for all}\; j \in \mathbb Z;$

\parindent=0mm \vspace{.1in}
(b)\quad $\bigcup_{j\in \mathbb Z}V_j\;\text{is dense in}\;L^2({\mathbb R});$

\parindent=0mm \vspace{.1in}
(c)\quad $\bigcap_{j\in \mathbb Z}V_j=\{0\};$

\parindent=0mm \vspace{.1in}
(d)\quad $f(x) \in V_j\; \text{if and only if}\;f(2Nx) \in V_{j+1}\; \text{for all}\; j \in \mathbb Z;$

\parindent=0mm\vspace{.1in}

(e)~ There exists a function $\phi$ in $V_0$ such that $\left\{ \phi (x- \lambda ): \lambda \in \Lambda\right\}$, is a complete orthonormal basis for $V_0$.

\parindent=0mm \vspace{.1in}
It is worth noticing that, when $N = 1$, one recovers from the definition above the  standard definition of one dimensional  multiresolution analysis with dyadic dilation. When, $N > 1$, the dilation factor of  $2N$ ensures that $2N\Lambda\subset {\mathbb Z} \subset \Lambda$.

\parindent=8mm \vspace{.1in}
For every $j\in\mathbb Z$, define $W_{j}$ to be the orthogonal complement of $V_{j}$ in $V_{j+1}$. Then we have
$$V_{j+1}=V_{j}\oplus W_{j}\quad \text{and}\quad W_{\ell}\perp W_{\ell^\prime}\quad \text{if}~\ell\ne \ell^\prime.\eqno(2.4) $$

\parindent=0mm \vspace{.1in}
It follows that for $j>J$,
$$V_{j}=V_{J}\oplus \bigoplus_{\ell=0}^{j-J-1}W_{j-\ell}\,,\eqno(2.5) $$

\parindent=0mm \vspace{.1in}
where all these subspaces are orthogonal. By virtue of condition (b) in the Definition 2.2, this implies

$$L^2({\mathbb R})=\bigoplus_{j\in\mathbb Z}W_{j},\eqno(2.6) $$

\parindent=0mm \vspace{.1in}
a decomposition of $L^2({\mathbb R})$ into mutually orthogonal subspaces.

\parindent=8mm \vspace{.1in}
As in the standard case, one expects the existence of $2N -1$ number of functions so that their translation by elements of $\Lambda$ and dilations by the integral powers of $2N$ form an orthonormal basis for $L^2({\mathbb R})$.

\parindent=0mm \vspace{.1in}
{\bf{Definition 2.3.}} A set of functions $\left\{\psi_{1}, \psi_{1},\dots,\psi_{2N-1}\right\}$ in $L^2({\mathbb R})$ is said to be a {\it set of  basic wavelets} associated with the nonuniform multiresolution analysis $\left\{V_{j}: j\in\mathbb Z\right\}$ if the family of functions $\left\{ \psi_{\ell}(x-\lambda):1\le \ell\le 2N-1, \lambda\in \Lambda\right\}$ forms an orthonormal basis for $W_{0}$.

\parindent=0mm \vspace{.2in}
{\bf{3.  Riesz Bases of Translates}}

\parindent=0mm \vspace{.1in}
{\bf{Lemma 3.1.}} {\it Let $\phi, \tilde \phi \in L^2({\mathbb R})$ be given. Then $\big\{\phi(x-\lambda):\lambda \in \Lambda\big\}$ is biorthogonal to $\big\{\tilde \phi(x-\lambda):\lambda \in \Lambda\big\}$ if and only if}

$$\sum_{\lambda \in \Lambda}\hat\phi(\xi+\lambda)\overline{\hat{\tilde \phi}(\xi+\lambda)}=1\quad a.e\; \xi \in \mathbb R.$$

\parindent=0mm \vspace{.0in}
{\it Proof.} For $\lambda, \sigma \in \Lambda$, it follows that $\left\langle \phi(x-\lambda), \tilde \phi(x-\sigma)\right\rangle=\delta_{\lambda, \sigma} \Leftrightarrow\left\langle \phi, \tilde \phi(x-\sigma)\right\rangle=\delta_{0, \sigma}$. Moreover, we have
$$\begin{array}{rcl}
\big\langle \phi, \tilde \phi(x-\sigma)\big\rangle&=&\left\langle \hat \phi,\hat{ \tilde \phi}(x-\sigma)\right\rangle\\
&=&\displaystyle\int_{\mathbb R}\hat\phi(\xi)\overline{\hat{\tilde \phi}(\xi)}e^{-2\pi i\sigma \xi}d\xi\\
&=&\displaystyle\int_0^{1/2}\left\{\sum_{p \in \mathbb Z}\hat\phi\left(\xi+\frac{p}{2}\right)\overline{\hat{\tilde \phi}\left(\xi+\frac{p}{2}\right)}e^{\pi i\sigma p}\right\}e^{-2\pi i\sigma \xi}d\xi.
\end{array}$$

\parindent=0mm \vspace{.1in}
Using the fact that $\big\{e^{-2\pi i\sigma \xi}:\sigma \in \Lambda\big\}$ is an orthonormal basis of $L^2\left[0, \frac{1}{2}\right)$, we get the desired result. \quad \fbox

\parindent=8mm \vspace{.1in}

We now provide a sufficient condition for the translates of a function to be linearly independent.

\parindent=0mm \vspace{.2in}

{\bf{Lemma 3.2.}} {\it Let $\phi \in L^2({\mathbb R})$. Suppose there exists two constants $A, B>0$ such that
$$A \le\sum_{\lambda \in \Lambda}\left|\hat\phi(\xi+\lambda)\right|^2\le B\quad for\; a.e \; \xi \in {\mathbb R}.\eqno (3.1)$$

\parindent=0mm \vspace{.0in}
Then, the collection  $\big\{\phi(x-\lambda):\lambda \in \Lambda\big\}$ is linearly independent.}

\parindent=0mm \vspace{.1in}
{\it Proof.} For the proof of the lemma, it is sufficient to find another function say $\tilde \phi$ whose translates are biorthogonal to  $\phi$. To do this, we define the function $\tilde \phi$ by
$$\hat {\tilde \phi}(\xi)=\dfrac{\hat \phi(\xi)}{\displaystyle\sum_{\lambda \in \Lambda}\left|\hat\phi(\xi+\lambda)\right|^2}.$$

\parindent=0mm \vspace{.0in}
Equation (3.1) implies that  $\tilde \phi$ is well defined and

$$\begin{array}{rcl}
\displaystyle\sum_{\sigma \in \Lambda}\hat\phi(\xi+\sigma)\overline{\hat{\tilde \phi}(\xi+\sigma)}&=&\displaystyle\sum_{\sigma \in \Lambda}\hat\phi(\xi+\sigma)\dfrac{\overline{\hat \phi(\xi+\sigma)}}{\displaystyle\sum_{\lambda \in \Lambda}\left|\hat\phi(\xi+\lambda+\sigma)\right|^2}\\
&=&\dfrac{\displaystyle\sum_{\sigma \in \Lambda}\left|\hat\phi(\xi+\sigma)\right|^2}{\displaystyle\sum_{\nu \in \Lambda}\left|\hat\phi(\xi+\nu)\right|^2}\\
&=&1.
\end{array}$$

\parindent=0mm \vspace{.1in}
Applying Lemma 3.1, it follows that $\big\{\phi(x-\lambda):\lambda \in \Lambda\big\}$ is linearly independent. This completes the proof of the Lemma.\quad \fbox\\

\parindent=0mm \vspace{.1in}
{\bf{Lemma 3.3.}} {\it Suppose that the scaling function $\phi$ satisfies inequality (3.1). Let $f=\sum_{\lambda \in \Lambda}h_\lambda\phi(x-\lambda)$, where $f \in {\text{span}}\big\{\phi(x-\lambda):\lambda \in \Lambda\big\}$ and $\big\{h_\lambda\big\}$ is a finite sequence. Define the Fourier transform of $h$ by $\hat h(\xi)=\displaystyle\sum_{\lambda \in \Lambda}h_\lambda e^{-2\pi i\lambda \xi}$. Then}
$$A\int_0^{1/2}\big|\hat h(\xi)\big|^2d\xi\le\big\|f\big\|_2^2\le B\int_0^{1/2}\big|\hat h(\xi)\big|^2d\xi. $$

\parindent=0mm \vspace{.0in}
{\it Proof.} By using Placherel's  theorem, we obtain
$$\begin{array}{rcl}
\displaystyle\int_{\mathbb R}\big|f(x)\big|^2dx&=&\displaystyle\int_{\mathbb R}\left|\sum_{\lambda \in \Lambda}h_\lambda\phi(x-\lambda)\right|^2dx\\
&=&\displaystyle\int_{\mathbb R}\left|\sum_{\lambda \in \Lambda}h_\lambda \hat\phi(\xi)e^{-2\pi i\lambda \xi}\right|^2d\xi\\
&=&\displaystyle\int_{\mathbb R}\big|\hat\phi(\xi)\big|^2\left|\sum_{\lambda \in \Lambda}h_\lambda e^{-2\pi i\lambda \xi}\right|^2d\xi\\
&=&\displaystyle\int_{\mathbb R}\big|\hat\phi(\xi)\big|^2\big|\hat h(\xi)\big|^2d\xi\\
&=&\displaystyle\int_0^{1/2}\sum_{p \in \mathbb Z}\left|\hat\phi\left(\xi+\frac{p}{2}\right)\right|^2 \left|\hat h(\xi)\right|^2d\xi.
\end{array}$$

\parindent=0mm \vspace{.0in}
Using identity (3.1), we get the desired result.\quad \fbox

\parindent=0mm \vspace{.2in}

{\bf{Theorem 3.4.}} {\it Let $\big\{\phi(x-\lambda):\lambda \in \Lambda\big\}$ be a Riesz basis for its closed linear {\text{span}}. Assume that there exists a function $\big\{\tilde \phi(x-\lambda):\lambda \in \Lambda\big\}$ which is biorthogonal to $\big\{\phi(x-\lambda):\lambda \in \Lambda\big\}$. Then, for every $f \in \overline{{\text{span}}}\left\{\phi(x-\lambda):\lambda \in \Lambda\right\}$, we have
$$ f =\sum_{\lambda \in \Lambda}\left\langle f, \tilde \phi(x-\lambda)\right\rangle\phi(x-\lambda); \eqno (3.2)$$

\parindent=0mm \vspace{.0in}
and there exists constants $A, B>0$ such that }
$$A\big\|f\big\|_2^2\le \sum_{\lambda \in \Lambda}\left|\left\langle f, \hat{\tilde \phi}(\xi-\lambda)\right\rangle\right|^2\le B\big\|f\big\|_2^2.\eqno (3.3)$$

\parindent=0mm \vspace{.1in}
{\it Proof.} We first prove (3.2) and (3.3) for any $f \in {\text{span}}\big\{\phi(x-\lambda):\lambda \in \Lambda\big\}$ and then generalize it to $\overline{{\text{span}}}\big\{\phi(x-\lambda):\lambda \in \Lambda\big\}$. Assume that $f \in {\text{span}}\big\{\phi(x-\lambda):\lambda \in \Lambda\big\}$, then there exists a finite sequence $\big\{h_\lambda: \lambda \in \Lambda\big\}$ such that $f =\sum_{\lambda \in \Lambda}h_\lambda \phi(x-\lambda)$. Also, the biorthogonality  condition implies that
$$\begin{array}{rcl}
\big\langle f, \tilde \phi(x-\sigma)\big\rangle&=&\left\langle \displaystyle\sum_{\lambda \in \Lambda}h_\lambda \phi(x-\lambda), \tilde \phi(x-\sigma)\right\rangle\\
&=&\displaystyle\sum_{\lambda \in \Lambda}h_\lambda\big\langle  \phi(x-\lambda), \tilde \phi(x-\sigma)\big\rangle\\
&=&h_\lambda,
\end{array}$$

\parindent=0mm \vspace{.0in}
which proves (3.2). In order to prove (3.3), we make the use of Lemma 3.3 to get

$$B^{-1}\big\|f\big\|_2^2\le\int_0^{1/2}\big|\hat h(\xi)\big|^2d\xi\le A^{-1}\big\|f\big\|_2^2. $$

\parindent=0mm \vspace{.0in}
Using the Placherel formula for Fourier series and the fact that $h_\lambda =\left\langle f, \tilde \phi(x-\lambda)\right\rangle$, we obtain
$$\int_0^{1/2}\big|\hat h(\xi)\big|^2d\xi=\sum_{\lambda \in \Lambda}\big|h_\lambda\big|^2=\sum_{\lambda \in \Lambda}\left|\big\langle f, \tilde \phi(x-\lambda)\big\rangle\right|^2.$$

\parindent=0mm \vspace{.0in}
This proves (3.3). We now generalize the results to $\overline{{\text{{\text{span}}}}}\left\{\phi(x-\lambda):\lambda \in \Lambda\right\}$.  For $f \in \overline{{\text{{\text{span}}}}}\left\{\tilde \phi(x-\lambda):\lambda \in \Lambda\right\}$, there exists a sequence $\{f_m:m \in \mathbb Z\}$ in ${\text{{\text{span}}}}\left\{\tilde \phi(x-\lambda):\lambda \in \Lambda\right\}$ such that $\|f_m-f\|_2 \to 0$ as $m \to \infty.$ Thus, for each $\lambda \in \Lambda$, we have
$$\big\langle f_m, \tilde \phi(x-\lambda)\big\rangle\to \big\langle f, \tilde \phi(x-\lambda)\big\rangle\quad \text{as}\; m \to \infty.$$

\parindent=0mm \vspace{.0in}
Hence, the result holds for each $f_m$. Thus, we have
\begin{align*}
\sum_{\lambda \in \Lambda}\left|\big\langle f, \tilde \phi(x-\lambda)\big\rangle\right|^2&=\sum_{\lambda \in \Lambda}\lim_{m \to \infty}\left|\big\langle f_m, \tilde \phi(x-\lambda)\big\rangle\right|^2\\
&=\lim_{m \to \infty}\sum_{\lambda \in \Lambda}\left|\big\langle f_m, \tilde \phi(x-\lambda)\big\rangle\right|^2\\
&\le B\lim_{m \to \infty}\big\|f_m\big\|_2^2\\
&=B\big\|f\big\|_2^2.\tag{3.4}
\end{align*}

\parindent=0mm \vspace{.0in}
Moreover, we have
\begin{align*}
\left\{\displaystyle\sum_{\lambda \in \Lambda}\left|\big\langle f_m, \tilde \phi(x-\lambda)\big\rangle\right|^2\right\}^{1/2}&\le \left\{\displaystyle\sum_{\lambda \in \Lambda}\left|\big\langle f_m-f, \tilde \phi(x-\lambda)\big\rangle\right|^2\right\}^{1/2}+\left\{\displaystyle\sum_{\lambda \in \Lambda}\left|\big\langle f, \tilde \phi(x-\lambda)\big\rangle\right|^2\right\}^{1/2}.
\end{align*}

\parindent=0mm \vspace{.0in}
As the upper bound in (3.3) holds for $f_m-f$ and lower bound for each $f_m$, we infer that

$$A^{1/2}\big\|f\big\|_2 \le B^{1/2}\big\|f_m-f\big\|_2+\left(\sum_{\lambda \in \Lambda}\left|\big\langle f_m, \tilde \phi(x-\lambda)\big\rangle\right|^2\right)^{1/2},$$

\parindent=0mm \vspace{.0in}
from which we conclude that
$$A\big\|f\big\|_2^2 \le \sum_{\lambda \in \Lambda}\left|\left\langle f, \tilde \phi(x-\lambda)\right\rangle\right|^2.\eqno(3.5)$$

\parindent=0mm \vspace{.0in}
Combining (3.4) and (3.5), we obtain (3.3). Similarly, we can prove (3.2) for $f \in \overline{{\text{{\text{span}}}}}\left\{\phi(x-\lambda):\lambda \in \Lambda\right\}$ and the proof of the theorem is complete.\quad \fbox\\

\parindent=0mm \vspace{.2in}
{\bf{4.  Properties of Biorthogonal  Wavelets on the Spectrum}}

\parindent=0mm \vspace{.1in}
Let $\{V_j: j \in \mathbb Z\}$ and $\{\tilde V_j: j \in \mathbb Z\}$ be biorthogonal NUMRA's with scling functions $\phi$ and $\tilde \phi$. Then there exists integral periodic functions $m_0$ and $\tilde m_0$ such that $\hat \phi(\xi)=m_0\left(\xi/2N\right)\hat \phi\left(\xi/2N\right)$ and $\hat {\tilde \phi}(\xi)=\tilde m_0\left(\xi/2N\right)\hat {\tilde \phi}\left(\xi/2N\right)$. Suppose there exists integral periodic functions $m_\ell$ and $\tilde m_\ell, 1 \le \ell \le 2N-1$ such that

$$\mathcal M(\xi) \overline{\tilde{\mathcal M}(\xi)}=I,\eqno(4.1)$$

\parindent=0mm \vspace{.0in}
where
$$\mathcal M(\xi)=\left(\begin{array}{cccc}
m_0\left(\dfrac{\xi}{2N}\right)&m_0\left(\dfrac{\xi}{2N}+\dfrac{1}{4N}\right)&\dots &m_0\left(\dfrac{\xi}{2N}+\dfrac{2N-1}{4N}\right)\\
m_1\left(\dfrac{\xi}{2N}\right)&m_2\left(\dfrac{\xi}{2N}+\dfrac{1}{4N}\right)&\dots &m_2\left(\dfrac{\xi}{2N}+\dfrac{2N-1}{4N}\right)\\
\vdots&\vdots&\ddots&\vdots\\
m_{2N-1}\left(\dfrac{\xi}{2N}\right)&m_{2N-1}\left(\dfrac{\xi}{2N}+\dfrac{1}{4N}\right)&\dots &m_{2N-1}\left(\dfrac{\xi}{2N}+\dfrac{2N-1}{4N}\right)
\end{array} \right)$$
and
 $$\tilde {\mathcal M}(\xi)=\left(\begin{array}{cccc}
\tilde m_0\left(\dfrac{\xi}{2N}\right)&\tilde m_0\left(\dfrac{\xi}{2N}+\dfrac{1}{4N}\right)&\dots &\tilde m_0\left(\dfrac{\xi}{2N}+\dfrac{2N-1}{4N}\right)\\
\tilde m_1\left(\dfrac{\xi}{2N}\right)&\tilde m_2\left(\dfrac{\xi}{2N}+\dfrac{1}{4N}\right)&\dots &\tilde m_2\left(\dfrac{\xi}{2N}+\dfrac{2N-1}{4N}\right)\\
\vdots&\vdots&\ddots&\vdots\\
\tilde m_{2N-1}\left(\dfrac{\xi}{2N}\right)&\tilde m_{2N-1}\left(\dfrac{\xi}{2N}+\dfrac{1}{4N}\right)&\dots &\tilde m_{2N-1}\left(\dfrac{\xi}{2N}+\dfrac{2N-1}{4N}\right)
\end{array} \right).$$

\parindent=0mm \vspace{.0in}
For $1 \le \ell \le 2N-1$, define the associated biorthgonal wavelets as $\psi_\ell$ and $\tilde \psi_\ell$ by
$$\hat \psi_\ell(\xi)=m_\ell\left(\xi/2N\right)\hat \phi\left(\xi/2N\right)\quad \text{and} \quad \hat {\tilde \psi}_\ell(\xi)=\tilde m_\ell\left(\xi/2N\right)\hat {\tilde \phi}\left(\xi/2N\right).$$

\parindent=0mm \vspace{.1in}
{\bf{Definition 4.1.}} A pair of NUMRA's $\{V_j: j \in \mathbb Z\}$ and $\{\tilde V_j: j \in \mathbb Z\}$ with scaling functions $\phi$ and $\tilde \phi$, respectively are said to be biorthogonal to each other if $\{\phi(\cdot-\lambda): \lambda \in \Lambda\}$ and $\{\tilde \phi(\cdot-\lambda): \lambda \in \Lambda\}$ are biorthogonal.

\parindent=0mm \vspace{.1in}
{\bf{Definition 4.2.}} Let $\phi$ and $\tilde \phi$ be scaling functions for biorthogonal NUMRA's. For each $j \in \mathbb Z$, define the operators $P_j$ and $\tilde P_j$ on $L^2(\mathbb R)$ by

$$P_j f=\sum_{\lambda \in \Lambda}\big \langle f, \tilde \phi_{j, \lambda}\big\rangle \phi_{j, \lambda}\quad \text{and} \quad \tilde P_j f=\sum_{\lambda \in \Lambda}\big \langle f, \phi_{j, \lambda}\big\rangle\tilde  \phi_{j, \lambda},$$

\parindent=0mm \vspace{.0in}
respectively. It is easy to verify that these operators are uniformly bounded on $L^2(\mathbb R)$ and both the series are convergent in  $L^2(\mathbb R)$.

\parindent=0mm \vspace{.1in}
{\bf{Remark 4.3.}} The operators $P_j$ and $\tilde P_j$ satisfy the following properties.

\parindent=0mm \vspace{.2in}
(a)\quad  $P_j f=f$ if and only if $f \in V_j$ and $\tilde P_j f=f$ if and only if $f \in \tilde V_j$.

\parindent=0mm \vspace{.2in}
(b)\quad  $\displaystyle\lim_{j \to \infty}\big\|P_j f-f\big\|_2=0$ and $\displaystyle\lim_{j \to -\infty}\big\|P_j f\big\|_2=0$ for every $f \in L^2(\mathbb R).$

\parindent=0mm \vspace{.2in}
{\bf{Theorem 4.4.}} {\it Let $\phi$ and $\tilde \phi$ be the scaling functions for biorthogonal NUMRA's and $\psi_\ell$ and $\tilde \psi_\ell, 1 \le \ell \le 2N-1$ be the associated wavelets satisfying (4.1). Then, we have the following}

\parindent=0mm \vspace{.1in}
(a)\quad  {\it $\big\{\psi_{\ell,0,\lambda}: \lambda \in \Lambda\big\}$ is biorthogonal to $\big\{\tilde \psi_{\ell,0,\sigma}:\sigma \in \Lambda\big\}$},

\parindent=0mm \vspace{.1in}
(b) \quad  $\big\langle \psi_{\ell,0,\lambda}, \tilde \phi_{0, \sigma}\big\rangle =\big\langle \tilde \psi_{\ell,0,\lambda}, \phi_{0, \sigma}\big\rangle,\quad {\text{for all}} \;\lambda, \sigma \in \Lambda$.

\parindent=0mm \vspace{.2in}
{\it Proof.} we have
\begin{align*}
&\sum_{p \in \mathbb Z}\hat \psi_\ell\left( \xi +\dfrac{p}{2}\right)\overline{\hat{\tilde \psi}_\ell\left( \xi +\dfrac{p}{2}\right)}\\
&=\displaystyle\sum_{p \in \mathbb Z}\left\{ m_\ell\left(\dfrac{\xi}{2N} +\dfrac{p}{4N}\right)\hat \phi\left(\dfrac{\xi}{2N} +\dfrac{p}{4N}\right) \overline{\tilde m_\ell\left(\dfrac{\xi}{2N} +\dfrac{p}{4N}\right)}\,\overline{\hat {\tilde \phi}\left(\dfrac{\xi}{2N} +\dfrac{p}{4N}\right)}\right\}\\
&=\sum_{s=0}^{2N-1}\sum_{p \in \mathbb Z}\left\{ m_\ell\left(\dfrac{\xi}{2N} +\dfrac{p}{2}+\dfrac{s}{4N}\right)\hat \phi\left(\dfrac{\xi}{2N} +\dfrac{p}{2}+\dfrac{s}{4N}\right) \overline{\tilde m_\ell\left(\dfrac{\xi}{2N} +\dfrac{p}{2}+\dfrac{s}{4N}\right)}\,\overline{\hat {\tilde \phi}\left(\dfrac{\xi}{2N} +\dfrac{p}{2}+\dfrac{s}{4N}\right)}\right\}\\
&=\sum_{s=0}^{2N-1}\left\{ m_\ell\left(\dfrac{\xi}{2N} +\dfrac{s}{4N}\right) \overline{\tilde m_\ell\left(\dfrac{\xi}{2N} +\dfrac{s}{4N}\right)}\right\}\\
&=1.
\end{align*}
Hence, by Lemma 3.1, $\big\{\psi_{\ell,0,\lambda}: \lambda \in \Lambda\big\}$ is biorthogonal to $\big\{\tilde \psi_{\ell,0,\lambda}:\lambda \in \Lambda\big\}$. This proves part (a). To prove part (b), we have for, $\lambda, \sigma \in \Lambda$
$$\begin{array}{rcl}
 \displaystyle\big\langle \psi_{\ell,0,\lambda}, \tilde \phi_{0, \sigma}\big\rangle&=&\big\langle \psi_\ell(x-\lambda), \tilde \phi(x-\sigma)\big\rangle\\\\
&=&\displaystyle\left\langle \hat\psi_\ell\,e^{-2\pi i\lambda}, \hat {\tilde \phi}\,e^{-2\pi i\sigma}\right\rangle\\
 &=&\displaystyle\int_{\mathbb R} m_\ell\left(\dfrac{\xi}{2N}\right)\hat \phi\left(\dfrac{\xi}{2N}\right)e^{-2\pi i\lambda\xi} \, \overline{\tilde m_0\left(\dfrac{\xi}{2N}\right)}\,\overline{\hat {\tilde \phi}\left(\dfrac{\xi}{2N}\right)}e^{2\pi i\sigma\xi} d\xi\\
&=&\displaystyle\int_0^{1/2}\sum_{p \in \mathbb Z}\left\{ m_\ell\left(\dfrac{\xi}{2N} +\dfrac{p}{4N}\right)\hat \phi\left(\dfrac{\xi}{2N} +\dfrac{p}{4N}\right) \right.\\
 &&\left.\qquad\qquad\qquad\qquad\times\overline{\tilde m_0\left(\dfrac{\xi}{2N} +\dfrac{p}{4N}\right)}\,\overline{\hat {\tilde \phi}\left(\dfrac{\xi}{2N} +\dfrac{p}{4N}\right)}\right\}e^{2\pi i (\sigma-\lambda)} d\xi\\
\end{array}$$

$$\begin{array}{rcl}
&=&\displaystyle\int_0^{1/2}\sum_{s=0}^{2N-1}\sum_{p \in \mathbb Z}\left\{ m_\ell\left(\dfrac{\xi}{2N} +\dfrac{p}{2}+\dfrac{s}{4N}\right)\hat \phi\left(\dfrac{\xi}{2N} +\dfrac{p}{2}+\dfrac{s}{4N}\right) \right.\\
 &&\left.\qquad\qquad\quad\times\overline{\tilde m_0\left(\dfrac{\xi}{2N} +\dfrac{p}{2}+\dfrac{s}{4N}\right)}\,\overline{\hat {\tilde \phi}\left(\dfrac{\xi}{2N} +\dfrac{p}{2}+\dfrac{s}{4N}\right)}\right\}e^{2\pi i (\sigma-\lambda)} d\xi\\
&=&\displaystyle\int_0^{1/2}\sum_{s=0}^{2N-1}\left\{ m_\ell\left(\dfrac{\xi}{2N} +\dfrac{s}{4N}\right) \overline{\tilde m_0\left(\dfrac{\xi}{2N} +\dfrac{s}{4N}\right)}\right\}e^{2\pi i (\sigma-\lambda)}d\xi\\
&=&0.
\end{array}$$

\parindent=0mm \vspace{.0in}
The dual one can also be shown equal to zero in a similar manner. This proves part (b) and hence completes the proof of the theorem. \quad \fbox\\

\parindent=0mm \vspace{.2in}
{\bf{Theorem 4.5.}} {\it Let $\phi,\tilde \phi,\psi_\ell$ and $\tilde \psi_\ell, 1 \le \ell \le 2N-1$ be as in Theorem 4.1. Let $\psi_0=\phi$ and  $\tilde \psi_0=\tilde \phi$. Then, for every $f \in L^2({\mathbb R})$, we have
\begin{align*}
P_1f=P_0f+\sum_{\ell=1}^{2N-1}\sum_{\lambda \in \Lambda}\big\langle f, \tilde \psi_{\ell, 0, \lambda}\big\rangle\psi_{\ell, 0, \lambda}\tag{4.2}
\end{align*}
and
\begin{align*}
\tilde P_1f=\tilde P_0f+\sum_{\ell=1}^{2N-1}\sum_{\lambda \in \Lambda}\big\langle f, \psi_{\ell, 0, \lambda}\big\rangle\tilde \psi_{\ell, 0, \lambda}.\tag{4.3}
\end{align*}

\parindent=0mm \vspace{.0in}
where the series in (4.2) and (4.3) converges in $L^2({\mathbb R})$.}

\parindent=0mm \vspace{.2in}
{\it Proof.} For the sake of convenience, we will only prove (4.2), as (4.3) is an easy consequence. In particular, we will prove it in the weak sense only. For this, let $f, g \in L^2({\mathbb R})$. Then, we have
$$\begin{array}{lcr}
\displaystyle\sum_{\ell=0}^{2N-1}\sum_{\lambda \in \Lambda}\left\langle f, \tilde \psi_{\ell, 0, \lambda}\right\rangle\overline{\big\langle g, \psi_{\ell, 0, \lambda}\big\rangle}&&\\
\qquad=\displaystyle\sum_{\ell=0}^{2N-1}\sum_{\lambda \in \Lambda}\left\{\int_{\mathbb R} \hat f(\xi)\overline{\hat {\tilde \psi}_\ell(\xi) }e^{2\pi i\lambda \xi}d\xi\right\}\left\{\int_{\mathbb R} \overline{\hat g(\xi)}\hat \psi_\ell(\xi) e^{-2\pi i\lambda \xi}d\xi\right\}&&\\
\qquad=\displaystyle\sum_{\ell=0}^{2N-1}\sum_{\lambda \in \Lambda}\left\{\int_0^{1/2}\sum_{p \in \mathbb Z} \hat f\left(\xi+\dfrac{p}{2}\right)\overline{\hat {\tilde \psi}_\ell\left(\xi+\dfrac{p}{2}\right) }e^{2\pi i\lambda \xi}d\xi\right\}&&\\
\qquad\qquad\qquad\qquad\qquad\qquad\qquad\qquad\times\left\{\displaystyle\int_0^{1/2}\sum_{q \in \mathbb Z} \overline{\hat g\left(\xi+\dfrac{q}{2}\right)}\hat \psi_\ell\left(\xi+\dfrac{q}{2}\right) e^{-2\pi i\lambda \xi}d\xi\right\}&&\\
\qquad=\displaystyle\sum_{\ell=0}^{2N-1}\int_0^{1/2}\left\{\sum_{p \in \mathbb Z} \hat f\left(\xi+\dfrac{p}{2}\right)\overline{\hat {\tilde \psi}_\ell\left(\xi+\dfrac{p}{2}\right) }\right\} \left\{\sum_{q \in \mathbb Z}\overline{\hat g\left(\xi+\dfrac{q}{2}\right)}\hat \psi_\ell\left(\xi+\dfrac{q}{2}\right) \right\}d\xi&&\\
\qquad=\displaystyle\int_0^{1/2}\sum_{\ell=0}^{2N-1}\left\{\sum_{p \in \mathbb Z} \hat f\left(\xi+\dfrac{p}{2}\right)\overline{\tilde m_\ell\left(\dfrac{\xi}{2N} +\dfrac{p}{4N}\right)}\,\overline{\hat {\tilde \phi}\left(\dfrac{\xi}{2N} +\dfrac{p}{4N}\right)}\right.&&\\
\qquad\qquad\qquad\qquad\qquad \qquad\qquad\times \left.\displaystyle\sum_{q \in \mathbb Z}\overline{\hat g\left(\xi+\dfrac{q}{2}\right)}m_\ell\left(\dfrac{\xi}{2N} +\dfrac{q}{4N}\right)\hat \phi\left(\dfrac{\xi}{2N} +\dfrac{q}{4N}\right) \right\}d\xi&&\\
\quad=\displaystyle\int_0^{1/2}\sum_{\ell=0}^{2N-1}\left\{\sum_{r=0}^{2N-1}\sum_{p^\prime\in \mathbb Z} \hat f\left(\xi+\dfrac{p^\prime}{2}N+\dfrac{r}{2}\right)\overline{\tilde m_\ell\left(\dfrac{\xi}{2N} +\dfrac{r}{4N}+\dfrac{p^\prime}{2}\right)}\,\overline{\hat {\tilde \phi}\left(\dfrac{\xi}{2N} +\dfrac{r}{4N}+\dfrac{p^\prime}{2}\right)}\right.&&\\\\
\qquad\qquad \times \displaystyle\sum_{s=0}^{2N-1}\sum_{q^\prime\in \mathbb N_0}\overline{\hat g\left(\xi+\dfrac{q^\prime}{2}N+\dfrac{s}{2}\right)}m_\ell\left(\dfrac{\xi}{2N} +\dfrac{s}{4N}+\dfrac{q^\prime}{2}\right)\left.\hat \phi\left(\dfrac{\xi}{2N} +\dfrac{s}{4N}+\dfrac{q^\prime}{2}\right) \right\}d\xi&&\\
\quad=\displaystyle\int_0^{1/2}\sum_{r=0}^{2N-1}\sum_{p^\prime\in \mathbb N_0}\sum_{s=0}^{2N-1}\sum_{q^\prime\in \mathbb N_0}\left\{\sum_{\ell=0}^{2N-1} \overline{\tilde m_\ell\left(\dfrac{\xi}{2N} +\dfrac{r}{4N}\right)}m_\ell\left(\dfrac{\xi}{2N} +\dfrac{s}{4N}\right)\right\}&&\\
\qquad\times\hat f\left(\xi+\dfrac{p^\prime}{2}N+\dfrac{r}{2}\right)\overline{\hat {\tilde \phi}\left(\dfrac{\xi}{2N} +\dfrac{r}{4N}+\dfrac{p^\prime}{2}\right)}\,\overline{\hat g\left(\xi+\dfrac{q^\prime}{2}N+\dfrac{s}{2}\right)}\hat \phi\left(\dfrac{\xi}{2N} +\dfrac{s}{4N}+\dfrac{q^\prime}{2}\right) d\xi&&\\
\quad=\displaystyle\int_0^{1/2}\sum_{p^\prime\in \mathbb N_0}\sum_{q^\prime\in \mathbb N_0}\sum_{s=0}^{2N-1}\hat f\left(\xi+\dfrac{p^\prime}{2}N+\dfrac{s}{2}\right)\overline{\hat {\tilde \phi}\left(\dfrac{\xi}{2N} +\dfrac{s}{4N}+\dfrac{p^\prime}{2}\right)}&&\\
\end{array}$$
$$\begin{array}{lcr}
\qquad\qquad\qquad\qquad\qquad \qquad \qquad\qquad\qquad \times \overline{\hat g\left(\xi+\dfrac{q^\prime}{2}N+\dfrac{s}{2}\right)}\hat \phi\left(\dfrac{\xi}{2N} +\dfrac{s}{4N}+\dfrac{p^\prime}{2}\right) d\xi&&\\
\quad=\displaystyle\sum_{s=0}^{2N-1}\int_0^{{s+1}/2}\sum_{p^\prime\in \mathbb N_0}\sum_{q^\prime\in \mathbb N_0}\hat f\left(\xi+\dfrac{p^\prime}{2}N\right)\overline{\hat {\tilde \phi}\left(\dfrac{\xi}{2N} +\dfrac{p^\prime}{2}\right)}\,\overline{\hat g\left(\xi+\dfrac{q^\prime}{2}N\right)}\hat \phi\left(\dfrac{\xi}{2N} +\dfrac{p^\prime}{2}\right) d\xi. \,(4.5)
\end{array}$$

\parindent=0mm \vspace{.1in}
Furthermore, we have
$$\begin{array}{lcr}
\displaystyle\sum_{\lambda \in \Lambda}\left\langle f, \tilde \phi_{1, \lambda}\right\rangle\overline{\left\langle g, \phi_{1, \lambda}\right\rangle}&&\\
\qquad =\displaystyle\sum_{\lambda \in \Lambda}\left\{\int_{\mathbb R} \hat f(\xi)\overline{\hat {\tilde \phi}\left(\dfrac{\xi}{2N}\right) }e^{2\pi i \xi/2N}d\xi\right\}\left\{\int_{\mathbb R} \overline{\hat g(\xi)}\hat \phi\left(\dfrac{\xi}{2N}\right) e^{-2\pi i \xi/2N}d\xi\right\}&&\\
\qquad =\displaystyle\int_0^{1/2}\sum_{p\in \mathbb Z} \hat f\left(\xi+\dfrac{p}{2}N\right)\overline{\hat {\tilde \phi}\left(\dfrac{\xi}{2N}+\dfrac{p}{2}\right) }d\xi\int_0^{1/2}\sum_{q\in \mathbb Z} \overline{\hat g\left(\xi+\dfrac{q}{2}N\right)}\hat \phi\left(\dfrac{\xi}{2N}+\dfrac{q}{2}\right) d\xi\\
\qquad =\displaystyle\int_0^{1/2}\sum_{p\in \mathbb Z} \hat f\left(\xi+\dfrac{p}{2}N\right)\overline{\hat {\tilde \phi}\left(\dfrac{\xi}{2N}+\dfrac{p}{2}\right) }d\xi\int_0^{1/2}\sum_{q\in \mathbb Z} \overline{\hat g\left(\xi+\dfrac{q}{2}N\right)}\hat \phi\left(\dfrac{\xi}{2N}+\dfrac{q}{2}\right) d\xi\\
\qquad =\displaystyle\int_0^{1/2}\sum_{p\in \mathbb Z}\sum_{q\in \mathbb Z}  \hat f\left(\xi+\dfrac{p}{2}N\right)\overline{\hat {\tilde \phi}\left(\dfrac{\xi}{2N}+\dfrac{p}{2}\right) }\overline{\hat g\left(\xi+\dfrac{q}{2}N\right)}\hat \phi\left(\dfrac{\xi}{2N}+\dfrac{q}{2}\right) d\xi.\quad\qquad (4.5)
\end{array}$$

\parindent=0mm \vspace{.0in}
Combing (4.5) and (4.5), we get the desired result. \quad \fbox\\

\parindent=0mm \vspace{.2in}
{\bf{Theorem 4.6.}} {\it Let $\phi,\tilde \phi,\psi_\ell$ and $\tilde \psi_\ell, 1 \le \ell \le 2N-1$ be as in Theorem 4.1. Then, for every $f \in L^2({\mathbb R})$, we have
$$f=\sum_{\ell=1}^{2N-1}\sum_{j \in \mathbb Z}\sum_{\lambda \in \Lambda}\left\langle f, \tilde \psi_{\ell, j, \lambda}\right\rangle\psi_{\ell, j, \lambda}=\sum_{\ell=1}^{2N-1}\sum_{j \in \mathbb Z}\sum_{\lambda \in \Lambda}\big\langle f, \psi_{\ell, j, \lambda}\big\rangle\tilde \psi_{\ell, j, \lambda}, \eqno (4.6)$$

\parindent=0mm \vspace{.0in}
where the series converges in $L^2({\mathbb R})$.}

\parindent=0mm \vspace{.2in}
{\it Proof.} Using Remark 4.3 and Theorem 4.5, proof of Theorem 4.6 follows. \quad \fbox

\parindent=0mm \vspace{.1in}

{\bf{Theorem 4.7.}} {\it Let $\phi$ and $\tilde \phi$ be the scaling functions for biorthogonal NUMRA's and $\psi_\ell$ and $\tilde \psi_\ell, 1 \le \ell \le 2N-1$ be the associated wavelets satisfying the matrix condition (4.1). Then, the collection $\big\{\psi_{\ell, j, \lambda}: 1 \le \ell \le 2N-1, j \in \mathbb Z, \lambda \in \Lambda \big\}$ and $\big\{\tilde \psi_{\ell, j, \lambda}: 1 \le \ell \le 2N-1, j \in \mathbb Z, \lambda \in \Lambda \big\}$ are biorthogonal. Moreover, if
$$\big|\hat \phi(\xi)\big|\le C(1+|\xi|)^{-\frac{1}{2}-\epsilon},\,\big|\hat{\tilde \phi}(\xi)\big|\le C(1+|\xi|)^{-\frac{1}{2}-\epsilon},\,\big|\hat \psi_\ell(\xi)\big|\le C|\xi|\; and \;\big|\hat{\tilde \psi}(\xi)\big|\le C|\xi|,\eqno (4.7)$$

\parindent=0mm \vspace{.0in}
for some constant $C>0,\,\epsilon>0 $ and for a.e. $\xi \in \mathbb R$, then $\big\{\psi_{\ell, j, \lambda}: 1 \le \ell \le 2N-1, j \in \mathbb Z, \lambda \in \Lambda \big\}$ and $\big\{\tilde \psi_{\ell, j, \lambda}: 1 \le \ell \le 2N-1, j \in \mathbb Z, \lambda \in \Lambda \big\}$ form Riesz bases for $L^2({\mathbb R})$.}

\parindent=0mm \vspace{.2in}
{\it Proof.} First we show that $\big\{\psi_{\ell, j, \lambda}: 1 \le \ell \le 2N-1, j \in \mathbb Z, \lambda \in \Lambda \big\}$ and $\big\{\tilde \psi_{\ell, j, \lambda}: 1 \le \ell \le 2N-1, j \in \mathbb Z, \lambda \in \Lambda \big\}$ are biorthogonal to each other. For this, we will show that for each $\ell,\, 1\le \ell \le 2N-1$ and $j \in \mathbb Z$,
$$\big\langle \psi_{\ell, j, \lambda}, \tilde \psi_{\ell, j, \sigma}\big\rangle=\delta_{\lambda, \sigma}.\eqno(4.8)$$

\parindent=0mm \vspace{.1in}
In fact, we have already proved (4.8) for $j=0$. For $j\ne 0$, we have

$$ \big\langle \psi_{\ell, j, \lambda}, \tilde \psi_{\ell, j, \sigma}\big\rangle=\big\langle D_{-j} \psi_{\ell, 0, \lambda},  D_{-j}\tilde \psi_{\ell, 0, \sigma}\big\rangle=\big\langle\psi_{\ell, 0, \lambda}, \tilde \psi_{\ell, 0, \sigma}\big\rangle=\delta_{\lambda, \sigma}.$$

\parindent=0mm \vspace{.1in}
Also, for fixed $\lambda, \sigma \in \Lambda$  and $j, j^\prime \in \mathbb Z$ with $j<j^\prime$, we claim that
$$\big\langle \psi_{\ell, j, \lambda}, \tilde \psi_{\ell^\prime, j^\prime, \sigma}\big\rangle=0.\eqno(4.9)$$

\parindent=0mm \vspace{.0in}
As $\psi_{\ell, 0, \lambda} \in V_1$, hence $\psi_{\ell, j, \lambda}=D_{-j} \psi_{\ell, 0, \lambda}\in V_{j+1}\subseteq V_{j^\prime}$. Therefore, it is enough to show that $\tilde \psi_{\ell^\prime, j^\prime, \sigma}$ is orthogonal to every element of $V_{j^\prime}$. Let $f \in V_{j^\prime}$. Since $\big\{\phi_{j^\prime, \lambda}: \lambda \in \Lambda\big\}$ is a Riesz basis for $V_{j^\prime}$, hence there exists an $l^2$-sequence $\big\{d_\lambda: \lambda \in \Lambda\big\}$ such that $f=\sum_{\lambda \in \Lambda}d_\lambda\phi_{j^\prime, \lambda}$ in $L^2(\mathbb R)$. Using part (b) of Lemma 4.1, we have
$$\big\langle  \tilde \psi_{\ell^\prime, j^\prime, \sigma}, \phi_{j^\prime, \lambda}\big\rangle=\big\langle D_{-j^\prime} \tilde \psi_{\ell^\prime, 0, \sigma},  D_{-j^\prime}\phi_{0, \lambda}\big\rangle=0.$$

\parindent=0mm \vspace{.0in}
Therefore,
$$\big\langle  \tilde \psi_{\ell^\prime, j^\prime, \sigma}, f\big\rangle=\Big\langle \tilde \psi_{\ell^\prime, j^\prime, \sigma},  \displaystyle\sum_{\lambda \in \Lambda}d_\lambda\phi_{j^\prime, \lambda}\Big\rangle=\displaystyle\sum_{\lambda \in \Lambda}d_\lambda\big\langle \tilde \psi_{\ell^\prime, j^\prime, \sigma},  \phi_{j^\prime, \lambda}\big\rangle=0.$$

\parindent=0mm \vspace{.0in}
We now show that these two collections form Riesz bases for $L^2(\mathbb R)$. The linear independence is clear from the fact that these collections are biorthogonal to each other. So, we have to check the frame inequalities only i.e., there must exist constants $A, \tilde A, B, \tilde B>0$ such that
$$A\big\|f\big\|_2^2\le \sum_{\ell=1}^{2N-1}\sum_{j \in \mathbb Z}\sum_{\lambda \in \Lambda}\left|\big\langle f,  \psi_{\ell, j, \lambda}\big\rangle\right|^2\le B\big\|f\big\|_2^2,\quad \text{for all} \; f \in L^2(\mathbb R), \eqno (4.10)$$
and
$$\tilde A\big\|f\big\|_2^2\le \sum_{\ell=1}^{2N-1}\sum_{j \in \mathbb Z}\sum_{\lambda \in \Lambda}\left|\big\langle f, \tilde \psi_{\ell, j, \lambda}\big\rangle\right|^2\le \tilde B\big\|f\big\|_2^2,\quad \text{for all} \; f \in L^2(\mathbb R). \eqno (4.11)$$

\parindent=0mm \vspace{.1in}

Let us first check the existence of the upper bounds in (4.10) and (4.1 ). For this, we have
$$\begin{array}{rcl}
\displaystyle\sum_{\lambda \in \Lambda}\left|\big\langle f, \tilde \psi_{\ell, j, \lambda}\big\rangle\right|^2&=&\displaystyle\sum_{\lambda \in \Lambda}\left|\int_{\mathbb R} \hat f(\xi)(2N)^{-j/2}\,\overline{\hat \psi_\ell\big((2N)^{-j}\xi\big)}e^{2\pi i \lambda (2N)^{-j}\xi}d\xi\right|^2\\
&=&(2N)^{-j}\displaystyle\sum_{\lambda \in \Lambda}\left|\int_0^{1/2} \sum_{p \in \mathbb Z}\hat f\left(\xi+(2N)^j\dfrac{p}{2}\right)\overline{\hat \psi_\ell\left((2N)^{-j}\xi+\dfrac{p}{2}\right)}e^{2\pi i \lambda(2N)^{-j}\xi}d\xi\right|^2\\
&=&\displaystyle\int_0^{1/2}\left| \sum_{p \in \mathbb Z}\hat f\left(\xi+(2N)^j\dfrac{p}{2}\right)\overline{\hat \psi_\ell\left((2N)^{-j}\xi+\dfrac{p}{2}\right)}\right|^2d\xi\\
\end{array}$$
$$\begin{array}{rcl}
&=&\displaystyle\int_0^{1/2}\left\{\sum_{p \in \mathbb Z}\left| \hat f\left(\xi+(2N)^j\dfrac{p}{2}\right)\right|^2\left|\hat \psi_\ell\left((2N)^{-j}\xi+\dfrac{p}{2}\right)\right|^{2\delta}\right\}\\
&&\qquad \qquad\qquad\qquad \qquad\qquad\times\left\{\displaystyle\sum_{q \in \mathbb Z}\left|\hat \psi_\ell\left((2N)^{-j}\xi+\dfrac{q}{2}\right)\right|^{2(1-\delta)}\right\}d\xi\\
&=&\displaystyle\int_{\mathbb R}\left| \hat f(\xi)\right|^2\big|\hat \psi_\ell\left((2N)^{-j}\xi\right)\big|^{2\delta}\sum_{q \in \mathbb Z}\left|\hat \psi_\ell\left((2N)^{-j}\xi+\dfrac{q}{2}\right)\right|^{2(1-\delta)}d\xi.
\end{array}$$

\parindent=0mm \vspace{.0in}
By our assumption (4.7), we have  $|\hat \psi_\ell(\xi)|\le C\left(1+|(2N)^{-1}\xi|\right)^{-1/2-\epsilon}$ and therefore, it follows that $\sum_{q \in \mathbb Z}\big|\hat \psi_\ell\left((2N)^{-j}\xi+q/2\right)\big|^{2(1-\delta)}$ is uniformly bounded if $\delta<2\epsilon(1+2\epsilon)^{-1}$. Thus, there exists a constant $C>0$ such that
$$\begin{array}{rcl}
\displaystyle\sum_{\lambda \in \Lambda}\left|\big\langle f, \tilde \psi_{\ell, j, \lambda}\big\rangle\right|^2&\le& C\displaystyle\int_{\mathbb R}\big| \hat f(\xi)\big|^2\sum_{j \in \mathbb Z}\big|\hat \psi_\ell\left((2N)^{-j}\xi\right)\big|^{2\delta}d\xi\\
&\le& C\sup\left\{\displaystyle\sum_{j \in \mathbb Z}\big|\hat \psi_\ell\left((2N)^{-j}\xi\right)\big|^{2\delta}: \xi \in [1, 2N]\right\}\big\|f\big\|_2^2.
\end{array}$$

\parindent=0mm \vspace{.0in}
Also for $\xi \in [1, 2N]$, we have
$$\begin{array}{rcl}
\displaystyle\sum_{j=-\infty}^0\big|\hat \psi_\ell\left((2N)^{-j}\xi\right)\big|^{2\delta}&\le&\displaystyle\sum_{j=-\infty}^0\dfrac{C^{2\delta}}{\big(1+|(2N)^{j-1}\xi|\big)^{\delta(1+2\epsilon)}}\\
&\le&\displaystyle\sum_{j=-\infty}^0\dfrac{C^{2\delta}}{(2N)^{(j-1)\delta(1+2\epsilon)}}\\
&\le&\displaystyle C^{2\delta}\dfrac{q^{\delta(1+2\epsilon)}}{1-(2N)^{-\delta(1+2\epsilon)}}.
\end{array}$$

\parindent=0mm \vspace{.0in}
Furthermore, we have
$$\sum_{j=1}^\infty\big|\hat \psi_\ell\left((2N)^{-j}\xi\right)\big|^{2\delta}\le\sum_{j=1}^\infty\big(C(2N)^{-j}|\xi|\big)^{2\delta}\le C^{2\delta}\sum_{j=1}^\infty (2N)^{(-j+1)2\delta}=C^{2\delta}\dfrac{1}{1-(2N)^{-2\delta}},
$$

\parindent=0mm \vspace{.0in}
and hence, it follows that $\sup\big\{\sum_{j \in \mathbb Z}\big|\hat \psi_\ell\left((2N)^{-j}\xi\right)\big|^{2\delta}: \xi \in [1, 2N]\big\}$ is finite. Therefore, there exist $B>0$ such that of (4.10) holds. Similarly, we can show for dual one also. The existence of lower bounds for both the cases can be shown in similar fashion. Using Theorem 4.6, it follows that if $f \in L^2(\mathbb R)$, then (4.6) holds. Thus, we have

$$\begin{array}{rcl}
\big\|f\big\|_2^2&=&\big\langle f, f \big\rangle\\
&=&\left\langle \displaystyle\sum_{\ell=1}^{2N-1}\sum_{j \in \mathbb Z}\sum_{\lambda \in \Lambda}\big\langle f, \tilde\psi_{\ell, j, \lambda}\big\rangle \psi_{\ell, j, \lambda}, f \right\rangle\\
&=&\displaystyle\sum_{\ell=1}^{2N-1}\sum_{j \in \mathbb Z}\sum_{\lambda \in \Lambda} \big\langle f, \tilde \psi_{\ell, j, \lambda}\big\rangle\big\langle \psi_{\ell, j, \lambda}, f \big\rangle\\
&\le&\left(\displaystyle\sum_{\ell=1}^{2N-1}\sum_{j \in \mathbb Z}\sum_{\lambda \in \Lambda} \left|\big\langle f, \tilde \psi_{\ell, j, \lambda}\big\rangle\right|^2\right)^{1/2}\left(\displaystyle\sum_{\ell=1}^{2N-1}\sum_{j \in \mathbb Z}\sum_{\lambda \in \Lambda} \left|\big\langle f, \psi_{\ell, j, \lambda}\big\rangle\right|^2\right)^{1/2}\\
&\le&(\tilde B)^{1/2}\big\|f\big\|_2\left(\displaystyle\sum_{\ell=1}^{2N-1}\sum_{j \in \mathbb Z}\sum_{\lambda \in \Lambda} \left|\big\langle f, \psi_{\ell, j, \lambda}\big\rangle\right|^2\right)^{1/2}.
\end{array}$$

\parindent=0mm \vspace{.0in}
Hence,
$$\dfrac{1}{\tilde B}\,\big\|f\big\|_2^2\le \sum_{\ell=1}^{2N-1}\sum_{j \in \mathbb Z}\sum_{\lambda \in \Lambda} \left|\big\langle f, \psi_{\ell, j, \lambda}\big\rangle\right|^2.$$

\parindent=0mm \vspace{.0in}
The dual case can be proved in similar lines. This completes the proof. \quad \fbox


\parindent=0mm \vspace{.2in}


\begin{thebibliography}{11}

{\small {

\bibitem{BG} M. Bownik and G. Garrigos, Biorthogonal wavelets, MRA's and shift-invariant spaces, Studia Math. 160, 231-248, (2004).

\bibitem {CDF} A. Cohen, I. Daubechies and J. C. Feauveau, Biorthogonal bases of compactly supported wavelets, Commun. Pure Appl. Math. 45, 485-560 (1992).

\bibitem{CW}  C. K. Chui and  J. Z. Wang, On compactly supported spline wavelets and a duality principle, Trans. Amer.  Math. Soc. 330(2), 903-915 (1992).

\bibitem{GN} J. P. Gabardo and M. Z. Nashed, Nonuniform multiresolution analysis and spectral pairs,  J. Funct. Anal., 158,  209-241 (1998)

\bibitem{GY} J. P. Gabardo and X. Yu, Wavelets associated with nonuniform multiresolution analysis and one-dimensional spectral pairs.  J. Math. Anal. Appl.  323, 798-817 (2006)

\bibitem{LC} R. Long and D. Chen, Biorthogonal wavelet bases on $\mathbb R^d$, Appl. Comput. Harmon. Anal. 2, 230-242, (1995).

\bibitem{SA} F. A. Shah and Abdullah, Nonuniform multiresolution analysis on local fields of positive characteristic. Comp. Anal.  Opert. Theory. 9, 1589-1608 (2015)

\bibitem{SB} F. A. Shah and M. Y. Bhat, Vector-valued nonuniform multiresolution analysis on local fields.  Int. J.  Wavelets, Multiresolut. Inf. Process. 13, 4 Article Id: 1550029 (2015).

\bibitem{SM} N. K. Shukla and S. Mittal, Wavelets on the spectrum, Numer. Funct. Anal.  Optimizat. 35(4), 461–486, (2014).

   }}


\end{thebibliography}
\end{document}